\documentclass[12pt]{amsart}
\usepackage[all]{xy}
\usepackage[dvips]{graphicx}
\usepackage[dvips]{color}
\usepackage{mathrsfs}
\usepackage[centertags]{amsmath}
\usepackage{amsfonts}
\usepackage{amssymb}
\usepackage{amsthm}

\newtheorem{theorem}{Theorem}

\newcommand{\C}{\mathbb{C}}
\newcommand{\R}{\mathbb{R}}
\newcommand{\E}{\mathbb{E}}
\newcommand{\K}{\mathbb{K}}

\def\S{\mathbb{S}}

\def \P{{\rm I\kern -2.2pt P\hskip 1pt}}

\begin{document}
\title[Random Polynomials ]{A review of some recent results on Random Polynomials over $\R$ and over $\C$.}
\author{Diego Armentano}
\date{}
\address{Centro de Matem\'atica, Universidad de la Rep\'ublica. Montevideo, Uruguay}

 \begin{abstract}
This article is divided in two parts. 
In the first part we review some recent results concerning the expected number of real roots of random system of polynomial equations. 
In the second part we deal with a different problem, namely, 
the distribution of the roots of certain complex random polynomials.
 We discuss a recent result in this direction, which shows that the associated
points in the sphere (via the stereographic projection) are
surprisingly well-suited with respect to the minimal logarithmic energy on the sphere. 
\end{abstract}

\maketitle

\section{Introduction}

Let us consider a system of $m$ polynomial equations in $m$ unknowns over a field $\K$,
\begin{eqnarray}\label{def:pol}
f_i(x):=\sum_{\|j\|\leq d_i}a_{j}^{(i)}x^j\quad (i=1,\ldots,m).
\end{eqnarray}
 The notation in (\ref{def:pol}) is the
following: $x:=(x_1,\ldots,x_m)$ denotes a point in $\K^m$, $j:=(j_1,\ldots,j_m)$ a multi-index
of non-negative integers, $\|j\|=\sum_{h=1}^mj_h$, $\,x^j=x^{j_1}\cdots x^{j_m}$,
$\,a_j^{(i)}=a^{(i)}_{j_1,\ldots,j_m}$, and $\,d_i$ is the degree of the polynomial $f_i$.

We are interested in the solutions of the system of equations
\begin{equation}\label{sofe}
f_i(x)=0\quad (i=1,\ldots,m),
\end{equation}
lying in some subset $V$ of $\K^m$.
Throughout this review we are mainly concerned with the case $\K=\R$ or $\K=\C$.

If we choose at random the coefficients $\{a^{(i)}_j\}$, then the solution of the system (\ref{sofe}) becomes a random subset of $\K^m$. 
This is the main object of this review. 

In the first part of this paper we focus on the real case. 
The main problem we consider is that of understanding $N^f(V)$: the number of solutions lying in the Borel subset $V$ of $\R^m$.

In the second part we deal with a different problem: 
How are the roots of complex polynomials distributed?  

This article is organized as follows:
\\
In \textit{Section \ref{sec:ENR}} we start with some historical remarks on random polynomials. After that we move to the case of random systems of equations. We mention some recent results for centered Gaussian distributions. 
In \textit{Section \ref{sec:SA}} we consider the non-centered case, which has also been called ``smooth-analysis'' in the last years. That is, we start with a fixed (non-random) polynomial system, then we perturb it with a polynomial noise, and we ask what can be said about the number of roots of the perturbed system.
In \textit{Section \ref{sec:BB}} we review a result which computes the expected number of roots of a random system of polynomial equations expressed in a different basis, namely, the Bernstein basis.
Finally in \textit{Section \ref{sec:FP}} we focus on the complex case. We discuss a recent result concerning the distribution of points in the sphere associated with roots of random complex polynomials.

 \vspace{5pt}

 This review follows the talk given by the author in the colloquium which was held the inauguration of the Franco-Uruguayan Institute of Mathematics, in Punta del Este, Uruguay, on December 2009. 

\section{The Number of Real Roots of Random Polynomials}\label{sec:ENR}

The study of the expectation of the number of real roots of a random polynomial started in the thirties with the work of Block and Polya \cite{BlP}. Further investigations were made by Littlewood and Offord \cite{LO}. However, the first sharp result is due to M. Kac (see Kac\cite{kac1,kac2}), who gives the asymptotic value 
$$
\E \left(N^f(\R)\right) \approx \frac{2}{\pi} \log d,\quad\mbox{as}\quad d\to +\infty,
$$
when the coefficients of the degree $d$ univariate polynomial $f$ are Gaussian centered independent random variables $N(0,1)$ (see the book by Bharucha--Reid and Sambandham
\cite{BRS}). 

The first important result in the study of real roots of random system of polynomial equations is due to Shub and Smale \cite{ShSm93} in $1992$, where the authors  computed the expectation of $N^f(\R^m)$ when the coefficients are
Gaussian centered independent random variables having variances:
\begin{eqnarray}\label{covss}
\E\left[ (a_j^{(i)})^2  \right]=\frac{d_i!}{j_1!\cdots j_m!\,(d_i-\|j\|)!}.
\end{eqnarray}
Their result was
\begin{eqnarray}\label{shsm}
\E\left(N^f(\R^m)\right)=\sqrt{d_1\cdots d_m},
\end{eqnarray}
that is, the square root of the B\'ezout number associated to the system.
The proof is based on a double fibration manipulation of the co-area formula.
Some extensions of their work, including new results for one
polynomial in one variable, can be found in Edelman--Kostlan\cite{EK}. There are
also other extensions to multi-homogeneous systems in McLennan\cite{Mcl},
and, partially, to sparse systems in Rojas\cite{Ro} and Malajovich--Rojas\cite{MRo}. A
similar question for the number of critical points of real-valued
polynomial random functions has been considered in Dedieu--Malajovich\cite{DM}.

The probability law of the Shub--Smale model defined in (\ref{covss}) has the simplifying property of being invariant under the action of the orthogonal group in $\R^m$. 
In Kostlan\cite{kos} one can find the classification of all Gaussian probability distributions over the coefficients with this geometric invariant property.

In 2005, Aza\"is and Wschebor gave a new and deep insight to this problem. The key point is using the Rice formula for random Gaussian fields (cf. Aza\"is--Wschebor\cite{JMAMWRice}). This formula allows one to extend the Shub--Smale result to other probability distributions over the coefficients. 
A general formula for $\E(N^f(V))$ when the random functions
$f_i\,(i=1,\ldots,m)$ are stochastically independent and their law
is centered and invariant under the orthogonal group on $\R^m$ can
be found in Aza\"is--Wschebor\cite{JMAMW2005}. This includes the Shub--Smale formula
(\ref{shsm}) as a special case. 
Moreover, Rice formula appears to be the instrument to consider a major problem in the subject which is to find the asymptotic distribution of $N^f(V)$ (under some normalization). The only published results of which the author is aware
concern asymptotic variances as $m\to +\infty$. 
(See Wschebor\cite{Wreview} for a detailed description in this direction and a simpler proof of Shub--Smale result).

\subsection{Non-centered Systems}\label{sec:SA}

The aim of this section is to remove the hypothesis that the
coefficients have zero expectation.

One way to look at this problem is to start with a non-random system of equations (the ``signal'')
\begin{eqnarray}\label{parte:det}
P_i(x)=0 \quad (i=1,\ldots,m),
\end{eqnarray}
perturb it with a polynomial noise $X_i(x)$ $(i=1,\ldots, m)$,
that is, consider
\begin{eqnarray*}
P_i(x)+X_i(x)=0 \quad (i=1,\ldots,m),
\end{eqnarray*}
and ask what one can say about the number of roots of the new system, or, how much the noise modifies the number of roots of the deterministic part. (For short, we denote $N^f=N^f(\R^m)$).

Roughly speaking, we prove in \textit{Theorem \ref{teo:smooth}} that
if the relation signal over noise is neither too big nor too small,
in a sense that will be made precise later on, there exist
positive constants  $C,\,\theta$, where $0< \theta <1$, such that
\begin{eqnarray}\label{smallo}
 \E(N^{P+X})\leq C\, \theta^m \E(N^{X}).
 \end{eqnarray}

Inequality (\ref{smallo}) becomes of interest if the starting non-random system (\ref{parte:det}) has a
large number of roots, possibly infinite, and $m$ is large. In this situation, the effect of adding
polynomial noise is a reduction at a geometric rate of the expected number of roots, as compared to the
centered case in which all the $P_i$'s are identically zero.

For simplicity we assume that the polynomial noise $X$ has the Shub-Smale distribution. However, one should keep in mind that the result can be extended to other orthogonally invariant distributions (cf. Armentano--Wschebor\cite{arm}).

Before the statement of \textit{Theorem \ref{teo:smooth}} below, we need to introduce some additional notations.

In this simplified situation, one only needs hypotheses concerning the relation between the signal $P$ and the Shub-Smale
noise $X$, which roughly speaking should neither be too small nor too big.

Since $X$ has the Shub-Smale distribution, from (\ref{covss}) we get 
$$
\mbox{Var}(X_i(x))= (1+\|x\|^2)^{d_i},\quad \forall x\in\R^m, \qquad (i=1,\ldots, m).
$$

Define
 \begin{eqnarray*}
 &&H(P_i):=\sup_{x\in\R^m}\left\{ (1+\|x\|)\cdot\left\| \nabla\left( \frac{P_i}{(1+\|x\|^2)^{d_i/2}}
 \right)(x) \right\|  \right\} ,  \\
 &&K(P_i):=\sup_{x\in\R^m\setminus \{0\} }\left\{ (1+\|x\|^2)\cdot\left| \frac{\partial}{\partial \rho}
 \left( \frac{P_i}{(1+\|x\|^2)^{d_i/2}} \right)(x) \right|  \right\} ,
 \end{eqnarray*}
for $i=1,\ldots,m$, where $\|\cdot\|$ is the Euclidean norm, and $\frac{\partial}{\partial \rho}$ denotes the derivative in the direction defined by $\frac{x}{\|x\|}$,
at each point $x\neq 0$.

For $r>0$, put:
\begin{eqnarray*}
L(P_i,r):=\inf_{\|x\|\geq r}\frac{P_i(x)^2}{(1+\|x\|^2)^{d_i}}\quad (i=1,\ldots,m) .
\end{eqnarray*}
One can check by means of elementary computations that for each $P$ as above, one has
$$
H(P)<\infty,\; K(P)<\infty.
$$
With these notations, we introduce the following hypotheses on the systems as $m$ grows:
 \begin{itemize}
 \item[$H_1)$]
\begin{subequations}
\begin{align}
A_m&=\frac{1}{m}\cdot\sum_{i=1}^m \frac{H^2(P_i)}{i}=\mbox{o}(1) \quad\mbox{as }\,m\to +
\infty\label{hm}\\
 B_m&=\frac{1}{m}\cdot\sum_{i=1}^m \frac{K^2(P_i)}{i}=\mbox{o}(1) \quad\mbox{as }\,m\to +
 \infty.\label{km}
\end{align}
\end{subequations}
\item[$H_2)$]
There exist positive constants $r_0,\,\ell$ such that if $r\geq
r_0$:
$$
L(P_i,r)\geq \ell \quad \mbox{for all } \, i=1,\ldots, m.
$$
\end{itemize}
\begin{theorem}\label{teo:smooth}
Under the hypotheses $H_1)$ and $H_2)$, one has
\begin{eqnarray}\label{Thsm}
\E(N^{P+X}) \leq C\, \theta^m \E(N^{X}),
\end{eqnarray}
where $C,\,\theta $ are positive constants, $0< \theta <1$.
\end{theorem}

\subsubsection{Remarks on the statement of \textit{Theorem \ref{teo:smooth}} }
\begin{itemize}
\item
It is obvious that our problem does not depend on the order in which the equations
$$
P_i(x)+X_i(x)=0 \quad (i=1,\ldots,m)
$$
appear. However, conditions (\ref{hm}) and (\ref{km}) in hypothesis $H_3)$ do depend on the order.
One can state them by saying that there exists an order $i=1,\ldots,m$ on the equations, such that (\ref{hm})
and (\ref{km}) hold true.
\item
Condition $H_1)$ can be interpreted as a bound on the quotient signal over noise. In fact, it concerns the
gradient of this quotient. In (\ref{km}) the radial derivative appears, which happens to decrease faster as
$\|x\|\to \infty$ than the other components of the gradient.

Clearly, if $H(P_i),\,K(P_i)$ are bounded by fixed
constants, (\ref{hm}) and (\ref{km}) are verified. Also, some of
them may grow as $m\to +\infty$ provided (\ref{hm}) and (\ref{km})
remain satisfied.
\item
Hypothesis $H_2)$ goes -- in some sense -- in the opposite direction:
For large values of $\|x\|$ we need a lower bound of the relation signal over noise.
\item
A result of the type of \textit{Theorem \ref{teo:smooth}} can not be
obtained without putting some restrictions on the relation signal
over noise. In fact, consider the system
\begin{eqnarray}\label{pertsigma}
P_i(x)+\sigma\,X_i(x)=0 \quad (i=1,\ldots,m),
\end{eqnarray}
where $\sigma$ is a positive real parameter.
If we let $\sigma\to
+\infty$, the relation signal over noise tends to zero and the
expected number of roots will tend to $\E(N^X)$.
On the other hand, if $\sigma \downarrow 0$, $\E(N^X)$ can have different behaviours. For example, if $P$ is a ``regular'' system´´, the expected value of the number of roots of
(\ref{pertsigma}) tends to the number of roots of
$P_i(x)=0,\;(i=1,\ldots,m)$, which may be much bigger than $\E(N^X)$. In this case, the relation signal over
noise tends to infinity.
\item As it was mentioned before we can extend \textit{Theorem \ref{teo:smooth}} to other orthogonally invariant distributions. However,  for the general version we need to add more hypotheses.
\end{itemize}

In the next paragraphs we are going to give two simple examples.

For the proof of \textit{Theorem \ref{teo:smooth}} and more examples with different noises see Armentano--Wschebor\cite{arm}.

\subsubsection{Some Examples}

We assume that the degrees $d_i$ are uniformly bounded.

For the first example, let
$$
P_i(x)=\|x\|^{d_i}- r^{d_i},
$$
where $d_i$ is even and $r$ is positive and remains bounded as $m$ varies.
Then, one has:
\begin{eqnarray*}
&&\frac{\partial}{\partial\rho} \left(\frac{P_i}{(1+\|x\|^2)^{d_i/2}}  \right)(x)  =   \frac{d_i\,
\|x\|^{d_i-1}+d_i\,r^{d_i}\,\|x\|}{ (1+\|x\|^2)^{\frac{d_i}{2} +1}}\leq \frac{d_i(1+r^{d_i})}
{(1+\|x\|^2)^{3/2}}\\
&&\nabla \left( \frac{P_i}{(1+\|x\|^2)^{d_i/2}}  \right) (x)  =  \frac{d_i\,\|x\|^{d_i-2}+d_i\,r^{d_i}}
{(1+\|x\|^2)^{\frac{d_i}{2}+1}}\, x
\end{eqnarray*}
which implies
\begin{eqnarray*}
\left\| \nabla \left( \frac{P_i}{(1+\|x\|^2)^{d_i/2}}  \right) (x)\right\|\leq \frac{d_i(1+r^{d_i})}
{(1+\|x\|^2)^{3/2}}.
\end{eqnarray*}
Again, since the degrees $d_1,\ldots,d_m$ are bounded by a constant
that does not depend on $m$, $H_1)$ follows. $H_2)$ also holds
under the same hypothesis.\\

Notice that an interest in this choice of the $P_{i}$'s lies in the fact
that obviously the system  $P_{i}(x)=0$ $(i=1,\ldots,m)$ has an infinite number of roots (all points in the sphere of
radius $r$ centered at the origin are solutions), but the expected number of roots of the perturbed system
is geometrically smaller than the Shub--Smale expectation, when $m$ is large.

Our second example is the following:
Let $T$ be a polynomial of degree $d$ in one variable that has $d$ distinct real roots. Define:
$$
P_i(x_1,\ldots,x_m)= T(x_i) \quad (i=1,\ldots,m).
$$
One can easily check that the system verifies our hypotheses, so that there exist $C,\,\theta$ positive
constants, $0<\theta <1$ such that
$$
\E(N^{P+X}) \leq C\,\theta^m d^{m/2},
$$
where we have used the Shub--Smale formula when the degrees are all
the same. On the other hand, it is clear that $N^P = d^m$ so that
the diminishing effect of the noise on the number of roots can be
observed. A number of variations of these examples for $P$ can be
constructed, but we will not pursue the subject here.

\subsection{Other Polynomial Basis}\label{sec:BB}

Up to now all probability measures were introduced in a particular basis, namely, the monomial basis $\{x^j \}_{\|j\|\leq d}$. However, in many situations, polynomial systems are expressed in different basis, for example, orthogonal polynomials, harmonic polynomials,Bernstein polynomials, etc.
So, it is a natural question to ask:
\textit{What can be said about $N^f(V)$ when the randomization is performed in a different basis?}

For the case of random orthogonal polynomials see Barucha-Reid and Sambandham\cite{BRS}, and Edelman--Kostlan\cite{EK} for  random harmonic polynomials.

In this section following Armentano--Dedieu\cite{AD} we give an answer to the  average number of real roots of a random system of equations expresed in the Bernstein basis. Let us be more precise:

The Bernstein basis is given by:
$$b_{d,k}(x) = \binom{d}{k} x^k (1 - x)^{d - k}, \ 0 \le k \le d,$$
in the case of univariate polynomials, and
$$b_{d,j}(x_1, \ldots , x_m) = \binom{d}{j} x_1^{j_1} \ldots x_m^{j_m}(1 - x_1 - \ldots - x_m)^{d - \| j \|}, \ \|j\| \le d,$$
for polynomials in $m$ variables, where $j = (j_1,  \ldots , j_m)$ is a multi-integer, and $\binom{d}{j}$ is the multinomial coefficient.

Let us consider the set of real polynomial systems in $m$ variables, 
$$
f_i(x_1, \ldots , x_m) = \sum_{\| j \| \leq d_i} a_j^{(i)} b_{d,j}(x_1, \ldots , x_m) \qquad (i=1,\ldots,m).
$$
Take the coefficients $a_j^{(i)}$ to be independent Gaussian standard random variables.

Define
$$\tau : \R^m \rightarrow {\mathbb{P}}\left(\R^{m+1}\right)$$
by
$$\tau (x_1, \ldots , x_m) = [x_1, \ldots , x_m, 1 - x_1 - \ldots - x_m].$$
Here ${\mathbb{P}}\left( \R^{m+1} \right)$ is the projective space associated with $\R^{m+1}$, $[y]$ is the class of the vector $y \in \R^{m+1}$, $y \ne 0$, for the equivalence relation defining this projective space. The (unique) orthogonally invariant probability measure in ${\mathbb{P}}\left( \R^{m+1}\right)$ is denoted by $\lambda_m$.

With the above notation the following theorem holds:
\begin{theorem} \label{teo:bernstein}
\begin{enumerate}
\item For any Borel set $V$ in $\R^m$ we have
$$\E \left(N^f(V)\right) =  \lambda_m (\tau (V)) \sqrt{d_1 \ldots d_m}.$$
In particular
\item $\E \left(N^f\right) =  \sqrt{d_1 \ldots d_m},$
\item $\E \left( N^f(\Delta^m) \right) =  {\sqrt{d_1 \ldots d_m}}/{2^m},$
where $$\Delta^m = \left\{ x \in \R^m \ : \ x_i \ge 0 \mbox{ and } x_1 + \ldots + x_m \le 1 \right\},$$
\item When $m=1$, for any interval $I = [\alpha, \beta] \subset \R$, one has
$$\E\left( N^f(I) \right) = \frac{\sqrt{d}}{\pi}\left(\arctan(2 \beta - 1) - \arctan(2 \alpha - 1)\right).$$
\end{enumerate}
\end{theorem}

The fourth assertion in \textit{Theorem \ref{teo:bernstein}} is deduced from the first assertion but it also can be derived from Crofton's formula (see for example Edelman--Kostlan\cite{EK}).

For the proof of \textit{Theorem \ref{teo:bernstein}} see Armentano--Dedieu\cite{AD}

\section{Distribution of Complex Roots of Random Polynomials}\label{sec:FP}

In this part we will see that points in the sphere associated with roots of Shub--Smale complex analogue random polynomials via the stereographic projection, are surprisingly well-suited with respect to the minimal logarithmic energy on the sphere.
That is, they provide a fairly good approximation to a classical minimization problem over the sphere, namely, the Elliptic Fekete points problem. 

 Next paragraphs follows closely Armentano--Beltr\'an--Shub\cite{ABS}, where one can find proofs and more detailed references.

 Given  $x_1,\ldots,x_N\in \S^2=\{x\in\R^3:\,\|x\|=1 \}$, let 
\begin{equation}\label{eq:fekete}
V(x_1,\ldots,x_N)=\ln\prod_{1\leq i<j\leq N}\frac{1}{\|x_i-x_j\|}=-\sum_{1\leq i<j\leq N}\ln\|x_i-x_j\|
\end{equation}
be the logarithmic energy of the $N$-tuple $x_1,\ldots,x_N$. 
Let
\[
V_N=\min_{x_1,\ldots,x_N\in\S^2}V(x_1,\ldots,x_N)
\]
denote the minimum of this function. 
$N$-tuples minimizing the quantity (\ref{eq:fekete}) are usually called Elliptic Fekete Points. 
The problem of finding (or even approximate) such optimal configurations is a classical problem (see White\cite{White1952} for its origins). 

During the last decades this problem has attracted much attention, and the number of papers concerning it has grown amazingly. The reader may see Kuijlaars-Saff\cite{KS} for a nice survey.

 In the list of Smale's problems for the XXI Century \cite{Sm2000}, problem number 7 reads:

\textit{Can one find $x_1,\ldots,x_N\in\S^2$ such that
\begin{equation}\label{eq:Smale7}
V(x_1,\ldots,x_N)-V_N\leq c\ln N,
\end{equation}
c a universal constant?}

More precisely, Smale demands a real number algorithm in the sense of Blum--Cucker--Shub--Smale\cite{BlCuShSm98} that with input $N$ returns a $N$-tuple $x_1,\ldots,x_N$ satisfying equation (\ref{eq:Smale7}), and such that the running time is polynomial on $N$.

One of the main difficulties when dealing with this problem is that the value of $V_N$ is not even known up to logarithmic precision. In Rakhmanov--Saff--Zhou\cite{RakhmanovSaffZhou1994} the authors proved that if one defines $C_N$ by 
\begin{equation}\label{VsubN}
V_N=-\frac{N^2}{4}\ln\left(\frac{4}{e}\right)-\frac{N\ln N}{4}+C_N N,
\end{equation}
then, 
\[
-0.112768770...\leq\liminf_{N\rightarrow\infty}C_N\leq\limsup_{N\rightarrow\infty}C_N\leq-0.0234973...
\]

Let $X_1,\ldots,X_N$ be independent random variables with common uniform distribution over the sphere.
One can easily show that the expected value of the function $V(X_1,\ldots,X_N)$ in this case is,
\begin{equation}\label{eq:uniform}
\mathbb{E}(V(X_1,\ldots,X_N))=-\frac{N^2}{4}\ln\left(\frac{4}{e}\right)+\frac{N}{4}\ln\left(\frac{4}{e}\right).
\end{equation}
Thus, this random choice of points in the sphere with independent uniform distribution already provides a reasonable approach to the minimal value $V_N$, accurate to the order of $O(N\ln N)$. 

On one side, this probability distribution has an important property, namely, invariance under the action of the orthogonal group on the sphere. However, on the other hand this probability distribution lacks on correlation between points. More precisely, in order to obtain well-suited configurations one needs some kind of repelling property between points, and in this direction independence is not favorable. 
Hence, it is a natural question whether other handy orthogonally invariant probability distributions may yield better expected values. 
Here is where complex random polynomials comes into account.

Given $z\in\C$, let 
$$
\hat z:= \frac{(z,1)}{1+|z|^2}\in\C\times\R \cong \R^3
$$
be the associated points in the Riemann Sphere, i.e. the sphere of radius $1/2$ centered at $(0,0,1/2)$. 
Finally, let
$$
X=2\hat{z}-(0,0,1)\in \S^2
$$
be the associated points in the unit sphere. 

Given a polynomial $f$ in one complex variable of degree $N$, we consider the mapping
\[
f\mapsto V(X_1,\ldots,X_N),
\]
where $X_i$ ($i=1,\ldots,N$) are the associated roots of $f$ in the unit sphere. Notice that this map is well defined in the sense that it does not depend on the way we choose to order the roots. 
\begin{theorem}\label{th:fekete}
Let $f(z)=\sum_{k=0}^N a_k z^k$ be a complex random polynomial, such that the coefficients $a_k$ are independent complex random variables, such that the real and imaginary parts of $a_k$ are independent (real) Gaussian random variables centered at $0$ with variance $\binom{N}{k}$. Then, with the notations above,
\[
{\mathbb E}\left(V(X_1,\ldots,X_N)\right)=-\frac{N^2}{4}\ln\left(\frac{4}{e}\right)-\frac{N\ln N}{4}+\frac{N}{4}\ln\frac{4}{e}.
\]
\end{theorem}
Comparing  \textit{Theorem \ref{th:fekete}} with equations (\ref{VsubN}) and (\ref{eq:uniform}), we see that the value of $V$ is surpringsingly small at points coming from the solution set of this random polynomials. 
More precisely, necessarily many random realizations of the coefficients will produce values of $V$ below the average and very close to $V_N$, possibly close enough to satisfy equation (\ref{eq:Smale7}).

Notice that, taking the homogeneous counterpart of $f$, \textit{Theorem \ref{th:fekete}} can be restated for random homogeneous polynomials and considering its complex projective solutions, under the identification of $\P(\C^2)$ with the Riemann sphere.
In this fashion, the induced probability distribution over the space of homogeneous polynomials in two complex variables  corresponds to the classical unitarily invariant Hermitian structure of the respective space 
(see Blum--Cucker--Shub--Smale\cite{BlCuShSm98}). Therefore, the probability distribution of the roots in $\P(\C^2)$ is invariant under the action of the unitary group. 

Is not difficult to prove that the unitary group action over $\P(\C^2)$ correspond to the special orthogonal group of the unit sphere. Hence, the distribution of the associated random roots on the sphere is orthogonally invariant. Thus, \textit{Theorem \ref{th:fekete}} is another geometric confirmation of the repelling property of the roots of this Gaussian random polynomials.

For a proof of \textit{Theorem \ref{th:fekete}} and a more detailed discussion on this account see Armentano--Beltr\'an--Shub\cite{ABS}. See also Shub--Smale\cite{ShSm93c}.

\providecommand{\bysame}{\leavevmode\hbox to3em{\hrulefill}\thinspace}
\providecommand{\MR}{\relax\ifhmode\unskip\space\fi MR }
\providecommand{\MRhref}[2]{%
  \href{http://www.ams.org/mathscinet-getitem?mr=#1}{#2}
}

\end{document}